\documentclass[fr,12pt]{article}
\usepackage[latin1]{inputenc}
\usepackage{amsmath,amsfonts,graphicx}

\def \R{\mathbb{R}} 
\def \N{\mathbb{N}}

\def \Z{\mathbb{Z}}
\def \N{\mathbb{N}}

\def \T{\mathbb{T}}

\def \AA{{\cal{A}}}
\def \AAb{\ov{\AA}}

\def \FF{{\cal{F}}}
\def \GG{{\cal{G}}}
\def \GGb{\ov{{\cal{G}}}}

\def \LL{{\cal{L}}}

\def \TT{{\cal{T}}}

\def \a{\alpha}
\def \b{\beta}
\def \g{\gamma}

\def \e{\varepsilon}
\def \g{\gamma}

\def \o{\omega}

\def \p{\phi}

\def \vp{\varphi}

\def \tS{{\widetilde{S}}}
\def \vv{\vert \vert}
\def \dun{\partial_1}
\def \dde{\partial_2}

\def \eqv{\Longleftrightarrow}

\def \ov{\overline}

\def \pr{{\rm pr}}
\def \Fb{\ov{F}}

\begin{document}




\bigskip
\bigskip
\centerline{\Large THE $C^0$ INTEGRABILIY OF SYMPLECTIC}
\centerline{\Large TWIST MAPS WITHOUT CONJUGATE POINTS}

\medskip
\centerline{Marc Arcostanzo}
\centerline{(Avignon University, LMA EA 2151, Avignon, France)} 

\bigskip
\bigskip

Let $d \geq 1$ be an integer, $\T^d$ the d-dimensional torus, and  
$$F : T^*\T^d \longrightarrow T^*\T^d$$
\noindent 
a twist map\footnote{this will be defined in part 1}. In this article, we assume that 
$F$ is without conjugate points$^1$ et we try to understand what consequences this might 
have on the dynamics of $F$.

\medskip
We first describe the periodic orbits of $F$. To state our result precisely, let  
$$ \ov{F} : T^*\R^d \longrightarrow T^*\R^d $$ 
be a lift of $F$ to $T^*\R^d$ (identified with $\R^d \times (\R^d)^*$). If $\o \in T^*\T^d$ 
is a periodic point of $F$ with period $N \in \N^*$ and $\ov{\o} = (x,p) 
\in \R^d \times (\R^d)^*$ a lift of $\o$, then for some $r \in \Z^d$ we have  
$$\ov{F}^N(x,p) = (x + r,p).$$ 
\noindent
Reciprocally, if this equality holds and $\o$ is the projection of $(x,p)$ on $T^*\T^d$, 
then $F^N(\o) = \o$, which means that $\o$ is a periodic point of $F$ with $N$ as a period. 
So we may see the following result as a way to describe the periodic orbits of $F$.

\bigskip
{\bf Theorem 1} : {\it Let $F$ be a twist map without conjugate points. For every $N \in \N^*$, 
for every $r \in \Z^d$, for every $x \in \R^d$, there is a unique $p \in (\R^d)^*$ 
such that $\ov{F}^N(x , p) = (x + r, p)$.} 

\bigskip
Let $x$ be a point on the torus $\T^d$. As a consequence of theorem 1, $F$ has a countable 
number of peridic orbits intersecting $T_x^*\T^d$. Each of them is determined by an integer 
$N \geq 1$ (which is a period of the orbit) and a vector $r \in \Z^d$ (we may call it the 
homotopy class of the orbit).

\medskip
We prove that if we fix $N$ an $r$ and let $x$ vary in $\T^d$, the set 
$$ \GGb_{N,r} = \{ (x,p) \in \T^d \times (\R^d)^* \ {\rm s.t.} \ \Fb^N(x,p) = (x+r,p) \} $$
\noindent
is a lift to $T^*\R^d$ of an invariant Lagrangian submanifold $\GG_{N,r}$ of $T^*\T^d$. So 
this gives rise to a sequence of Lagrangian submanifolds, each of them being a union of periodic 
orbits of $F$. It is natural to wonder if we can find other invariant Lagrangian submanifolds 
of $T^*\T^d$. It was suggested by M. Bialy in [Bi] that $F$ is without conjugate points if and 
only $T^*\T^d$ may be written as the union of $F$-invariant Lagrangian graphs. 

J. Cheng and Y. Sun showed that 
this holds true when $d=1$ : they proved in [Ch-Su] that $F$ is without conjugate points if and 
only if $T^*\T^1$ is foliated by continuous, closed, invariant curves that are not null-homotopic 
(a standard result due to Birkhoff states that these curves must be graphs over $\T^1$). Here we  
generalize their result in any dimension as follows :  

\bigskip
{\bf Theorem 2} : {\it Let $F : T^*\T^d \longrightarrow T^*\T^d$ be a twist map. Then $F$ is 
without conjugate points if and only if there is a continuous foliation of $T^*\T^d$ by 
Lispchitz, Lagrangian invariant graphs.}

\bigskip
The techniques used by Cheng et Sun do not carry over to the higher dimensional case. Our proof 
uses ideas coming from weak KAM and Aubry-Mather theory : each leaf of the foliation is a dual 
Aubry set associated to some cohomology class. This strategy was already used in [AABZ], where 
it is shown that a similar result holds for a class of Hamiltonian flows. However, many arguments 
used in the continuous setting have no analogue in the discrete case. For example, there is no 
useful numerical quantity (as the Hamiltonian in the continuous case) which is constant along the 
orbit of $F$.

\bigskip
This article is organized as follows : we recall some basic facts in section 1. The proof of theorem 
1 is given in section 2. In the three following sections, we show that the sets $\GG_{N,r}$ are 
invariant Lagrangian graphs of $T^*\T^d$ as well as dual Aubry sets associated to a cohomology class. 
Finally we give a proof of theorem 2 in the last two sections.   
 
\bigskip
\bigskip
{\Large 1. Twist maps without conjugate points}

\bigskip
\medskip
Let $d \geq 1$ be an integer. Denote by $\T^d = \R^d / \Z^d$ the d-dimensional torus. 
Let $\TT = \R^d \times \R^d$ and let $\TT^* = \R^d \times (\R^d)^*$ be the cotangent space  
of $\R^d$. Consider a generating function, that is a map $S : \TT \longrightarrow \R$ 
of class $C^2$ which satisfies the following two conditions 

\noindent
$(i) \ \forall r \in \Z^d, \forall (x,y) \in \TT, S(x+r,y+r) = S(x,y)$ ;

\noindent
$(ii)$ (`uniform twist condition', see [Bi-McK]) There is a real number $A > 0$ for which  

$$ \forall (x,y) \in \TT, \forall \xi \in \R^d, \ 
\sum _{i,j} {\partial^2S(x,y) \over \partial x_i \partial y_j}(x,y) \xi_i \xi_j \leq - A \vv \xi \vv^2.$$

\medskip
We may then define various notions of action. For example, 
the action of a finite sequence $\g = (x_0, x_1, \dots, x_n)$ with values in $\R^d$ is    
$$S(\g ) = S(x_0, x_1, \dots, x_n) = \sum_{k=0}^{n-1} S(x_k,x_{k+1}) 
.$$

\noindent
If we fix an integer $n \geq 2$ and two points $x_0$ and $x_n$ in $\R^d$, we can define 
the `action with fixed endpoints' as the map
$$S_{x_0,x_n,n} : (x_1, \dots x_{n-1}) \in (\R^d)^{n-1} \longmapsto  
S(x_0, x_1, \dots, x_{n-1}, x_n) \in \R.$$
\noindent
Its critical points are the finite sequences $(x_1, \dots x_{n-1})$ for which 
$$\forall k \in \{ 1, \dots, n-1 \}, \ \dde S(x_{k-1} , x_k) + \dun S(x_k , x_{k+1}) = 0.$$
\noindent
The sequence $(x_0, x_1, \dots, x_n)$ will be called a finite extremal sequence. An (infinite) 
sequence $(x_n)_{n \in \Z}$ with values in $\R^d$ is said to be extremal if it satisfies 
$$\forall n \in \Z, \ \dde S(x_{n-1} , x_n) + \dun S(x_n , x_{n+1}) = 0.$$

\medskip
Condition $(ii)$ implies (see [Go], chapter 4) that for every $x_0 \in \R^d$ and every 
$y_0 \in \R^d$, the maps  
$$ x \longmapsto \dde S(x , y_0) \ \ {\rm and} \ \  y \longmapsto \dun S(x_0 , y)$$
\noindent
are diffeomorphisms. As an immediate consequence, every finite extremal sequence may be uniquely 
extended to an infinite extremal sequence. In particular, for every $(x,y) \in \TT$, there is 
a unique extremal sequence $(x_n)_{n \in \Z}$ for which $x_0 = x$ and $x_1 = y$. We shall 
denote by 
$$\vp : (x,y) = (x_0, x_1) \in \TT \longmapsto  (x_1, x_2) \in \TT $$
\noindent 
the corresponding shift diffeomorphism.

\medskip
The generating function also gives rise to a symplectic diffeomorphism $F$ of $T^*\T^d$, the 
cotangent bundle of $\T^d$. Let $\Fb : \TT^* \longrightarrow \TT^*$ be the diffeomorphism 
(twist map) implicitely defined by  
$$ \Fb(x,p) = (x',p') \eqv p = -\dun S (x , x') \ {\rm and} \ p' = \dde S (x , x').$$
\noindent 
It turns out that $\Fb$ et $\vp$ are conjugated : the map  
$$\LL : (x,y) \in \TT \longmapsto (x, -\dun S(x,y)) \in \TT^*$$
\noindent
is a diffeomorphism for which $\Fb = \LL \circ \vp \circ \LL^{-1}$. The diffeomorphism $\Fb$ 
is exact symplectic, which means that $\Fb^* \a - \a = dS$, where 
$\a = \sum_{i=1}^d x_i dq_i$ is the Liouville 1-form on $T^*\R^d$. Note that condition $(i)$ 
implies that $\Fb$ is the lift to $\TT^*$ of a symplectic diffeomorphism $F$ of $T^*\T^d$. 
In this article, we are interested in the dynamics of $F$, and we will use $S$ as a useful tool 
for our study.

\medskip
As a matter of fact, condition $(ii)$ has strong consequences on the behaviour of $S$. For 
example, the following result may be shown (see [Go] page 105 or [McK-Me-St] page 568 for a proof).  

\medskip
{\bf Lemma 1.1} : {\it There exists $\a \in \R$, $\b \in \R$ and $\g > 0$ such that } 
$$ \forall (x,y) \in \TT, S(x,y) \geq \a + \b \vv x - y \vv + \g \vv x - y \vv^2.$$

\smallskip
As an immediate consequence, we can construct extremal sequences going trough two arbitrary 
points in $\R^d$.

\medskip
{\bf Lemma 1.2} : {\it For every $(x,y) \in \TT$, for every integer $N \geq 1$, there exists 
an extremal sequence $(x_n)_{n \in \Z}$ for which $x_0 = x$ and $x_N = y$}.

\medskip
{\it Proof} : we already know that this is the case when $N=1$. When $N \geq 2$, it suffices 
to show that the map $S_{x,y,N}$ has a critical point. In view of lemma 1.1, $S_{x,y,N}$ is 
coercive and therefore achieves its minimum at a point $(x_1, \dots, x_{N-1})$. We then 
extend the finite extremal sequence $(x_0 = x, x_1, \dots, x_{N-1}, x_N = y)$ to an (infinite) 
extremal sequence $(x_n)_{n \in \Z}$. 
\hfill{$\square$}

\medskip
Let us show that this extremal sequence is unique if we assume that $F$ is without conjugate 
points. Let $\pi : (x,p) \in T^*\T^d \longmapsto x \in \T^d$ be the canonical projection. 
For every $(x,p) \in T^*\T^d$, the vertical space at $(x,p)$ is  
$$ V(x,p ) = {\rm Ker} \left ( D \pi _{\vert T_{(x,p)}T^*\T^d} \right ) .$$

\medskip
{\bf Definition 1.3} : {\it $F$ is without conjugate points if  
$$ \forall (x,p) \in T^*\T^d, \forall n \in \Z^*, V(F^n(x,p)) \cap DF^n(x,p) \cdot V(x,p) \ 
= \ \{ 0 \}.$$}

\smallskip
{\bf Proposition 1.4} : {\it If $F$ is without conjugate points, then for every $(x,y) \in \TT$ 
and for every integer $N \geq 2$, the map $S_{x,y,N}$ has a unique critical point ; and at that 
point, $S_{x,y,N}$ achieves its minimum. }

\medskip
{\it Proof} : For the 'existence' part, we refer to the proof of lemma 1.2. Now assume by 
contradiction that  $S_{x,y,N}$ has (at least) two distinct critical points. It is shown in [Bi-McK] 
that if $F$ is without conjugate points, then every critical point of $S_{x,y,N}$ is in fact a 
strict local minimum. $S_{x,y,N}$ is then a coercive $C^2$ map with two distinct strict local minima. 
We can apply an existence theorem for saddle points in finite dimension (see [St], Theorem 1.1, 
page 74). It says that $S_{x,y,N}$ possesses a third critical point which is not a local minimum 
of $S_{x,y,N}$. This is a contradiction. \hfill{$\square$}

\medskip
{\bf Corollary 1.5} : {\it If $F$ is without conjugate points, then for every $(x,y) \in \TT$ 
and every integer $N \geq 1$, there is a unique extremal sequence $(x_n)_{n \in \Z}$ 
with $x_0 = x$ et $x_N = y$}. 

\medskip
{\bf Remark 1.6} : Assume that $F$ is without conjugate points. Let $(x_n)_{n \in \Z}$ be 
an extremal sequence, and let $k$ and $l$ be two integers with $l - k \geq 2$. It follows 
from proposition 1.4 that 
$$\forall (y_{k+1}, \dots , y_{l-1}) \in (\R^d)^{l-k-1}, \ 
S(x_k, x_{k+1}, \dots, x_{l-1}, x_l) \leq S(x_k, y_{k+1}, x_{k+1}, x_l).$$
\noindent
Equality holds if and only if $y_i = x_i$ for every $i \in \{ k+1, \dots, l-1 \}$. This 
means that an extremal sequence minimizes the action with fixed endpoints between any 
two of its points.

\bigskip
\bigskip
{\Large 2. Construction of periodic orbits}

\bigskip
\medskip
In this section, we prove theorem 1. Let us fix $r \in \Z^d$ and $N \in \N^*$. We have to show 
that for every $x \in \R^d$, there is a unique $p \in (\R^d)^*$ for which $\Fb^N(x , p) = (x + r, p)$. 
The change of variable $(x,y) = \LL^{-1}(x,p)$ and the relation $\Fb = \LL \circ \vp \circ \LL^{-1}$ 
lead to 
$$ \Fb^N(x , p) = (x + r, p) \Longleftrightarrow \vp^N(x,y) = (x+r,y+r).$$

\noindent
According to the results of the previous part, the only $y \in \R^d$ for which we may have 
$\vp^N(x,y) = (x+r,y+r)$ is $y = x_1$, where $(x_n)_{n \in \Z}$ is the unique extremal sequence 
satisfying $x_0 = x$ and $x_N = x+r$. And for this $y$, we have $\vp^N(x,y) = (x+r,y+r)$ if and 
only if si $x_{N+1} = x_1 + r$. It turns out that this equality holds, as it is a special case of 
the following general result.

\medskip
{\bf Proposition 2.1} : {\it Let $r \in \Z^d$ and $N \in \N^*$. If $(x_n)_{n \in \Z}$ is an 
extremal sequence such that $x_N = x_0 + r$, then $x_{n+N} = x_n + r$ for all $n \in \Z$.}  

\medskip
For the proof, we shall use a technique of metric geometry introduced by H. Busemann 
(see [Bu], section 32) when he was studying G-spaces without conjugate points.

\medskip
For every $(x,y) \in \TT$ and for every integer $N \in \N^*$, we denote by $\AA_N(x,y)$ 
the minimum of the function $S_{x,y,N}$ when $N \geq 2$, and $S(x,y)$ if $N=1$. As the minimum 
is attained at a single point, $\AA_N$ is a continuous function. We clearly have  
$\AA_N(x+r,y+r) = \AA_N(x,y)$ for every $(x,y) \in \TT$ and every $r \in \Z^d$. 

\medskip
{\bf Lemma 2.2} : {\it For every $x$, $y$, $z$ in $\R^d$, for every  $N$, $N'$ in $\N^*$, 
the following triangular inequality holds :   
$$\AA_{N+N'}(x,z) \leq \AA_N(x,y) + \AA_{N'}(y,z).$$
\noindent
Moreover, one has equality if and only if $y = w_N$, where $(w_n)$ is the extremal sequence 
for which $w_0 = x$ et $w_{N+N'} = z$}.

\medskip
{\it Proof} : Let $(x_n)$ be the extremal sequence with $x_0 = x$ and $x_N = y$, and  
$(y_n)$ the extremal sequence with $y_N = y$ and $y_{N+N'} = z$. So we have 
$$ \AA_N(x,y) = S(x_0, x_1, \dots, x_N) \ {\rm and} \ 
\AA_{N'}(y,z) = S(y_N, y_{N+1}, \dots, y_{N+N'}).$$

Let $(z_n)$ be the sequence defined by  
$$z_n = \begin{cases} x_n \ {\rm for} \ n \leq N \\
 y_n \ {\rm for} \ n \geq N \end{cases}$$
\noindent
As we have $z_0 = x_0 = x$ and $z_{N+N'} = y_{N+N'} = z$, the definition 
of $\AA_{N+N'}(x,z)$ implies that 
$$ \AA_{N+N'}(x,z) \leq S(z_0, z_1, \dots , z_{N+N'}) = S(x_0, x_1, \dots, x_N) + 
S(y_N, y_{N+1}, \dots, y_{N+N'}), $$
\noindent
whence the inequality $\AA_{N+N'}(x,z) \leq \AA_N(x,y) + \AA_{N'}(y,z)$. 

If equality holds, then $S_{x,z,N+N'}$ achieves its minimum at $(z_1, z_2, \dots, z_{N+N'-1})$. 
But $S_{x,z,N+N'}$ achieves its minimum at a unique point, namely $(w_1, w_2, \dots, w_{N+N'-1})$. 
So we must have $z_N = w_N$, and therefore $y = w_N$. \hfill{$\square$}

\medskip
Consider the function 
$$ f : x \in \R^d \longmapsto \AA_N(x,x+r) \in \R.$$

As $f$ is continuous and $\Z^d$-periodic, there exists two points $a$ and $b$ in $\R^d$ 
with $f(a) = \min_{\R^d}f$ and $f(b) = \max_{\R^d}f$. We first establish proposition 2.1 
for a particular extremal sequence.

\medskip
{\bf Lemma 2.3} : {\it The extremal sequence $(x_n)$ for which $x_0 = b$ and $x_N = x_0 + r$ 
satisfies  
$$\forall n \in \Z, \ x_{n+N} = x_n + r.$$}

\smallskip
{\it Proof} : Using the periodicity of $\AA_{2N}$ and the triangular inequality, we get   
$$\AA_{2N}(x_0,x_{2N}) = \AA_{2N}(x_0 + r,x_{2N} + r)
\leq \AA_{N}(x_0 + r,x_{2N}) + \AA_{N}(x_{2N},x_{2N} + r),$$
\noindent
so that
$$\AA_{2N}(x_0,x_{2N}) \leq \AA_{N}(x_N,x_{2N}) + f(x_{2N}).$$ 
\noindent
As extremal sequences are action-minimizing (see remark 1.6), we also have
$$\AA_{2N}(x_0,x_{2N}) = \AA_N(x_0,x_N) + \AA_N(x_N,x_{2N}) = f(b) + \AA_N(x_N,x_{2N}),$$ 
\noindent
so that the last inequality yields  
$$\AA_{2N}(x_0,x_{2N}) \leq \AA_{2N}(x_0,x_{2N}) - f(b) + f(x_{2N}) 
\leq \AA_{2N}(x_0,x_{2N}),$$
\noindent
because $f$ achieves its maximum at $b$. This implies that equality holds in all the previous 
inequalities. Lemma 2.2 then tells us that $x_{2N} = y_N$, where $(y_n)$ is the unique extremal 
sequence with $y_0 = x_0 + r$ and $y_{2N} = x_{2N} + r$. 

As the extremal sequences $(y_n)$ and $(x_n + r)$ are equal at $n = 0$ and $n = 2N$, 
corollary 1.5 implies that they are equal for all $n$. So we have $y_N = x_N + r$, and therefore  
$x_{2N} = y_N = x_N + r$. Now the two extremal sequences $(x_{n+N})$ and $(x_n + r)$ 
are equal at $n = 0$ and $n = N$, so they are equal. \hfill{$\square$} 

\medskip
{\bf Lemma 2.4} : {\it The fonction $f$ is constant.}

\medskip
{\it Proof} : We only need to show that $\max_{\R^d} f = f(b) \leq f(a) = \min_{\R^d} f$. 
From the preceding lemma, we have $x_{nN} = x_0 + nr = b + nr$ for all integer $n$, so that  
$$ \forall n \geq 1, \AA_{nN}(b,b + nr) = n \AA_N(b,b + r) = n f(b).$$
\noindent
On the other hand, the triangular inequality implies that for every $n \geq 3$, 
$$\AA_{nN}(b , b + nr) \leq \AA_N(b , a + r) + \sum_{i=1}^{n-2} \AA_N(a + ir , a + (i+1)r) 
+ \AA_N(a + (n-1)r , b + nr).$$
\noindent
These two relations and the fact that $\AA_N$ is $\Z^d$-invariant lead to  
$$ n f(b) \leq \AA_N(b , a + r) + (n-2) f(a) + \AA_N(a , b + r). $$
\noindent
When we divide by $n$ and let $n$ go to infinity, we obtain $f(b) \leq f(a)$. 
\hfill{$\square$}

\medskip
As the function $f$ achieves its maximum at every point, the conclusion of lemme 2.3 holds 
for every $b \in \R^d$. This ends the proof of proposition 2.1 and the proof of theorem 1.

\medskip
{\bf Corollary 2.5} : {\it If $F$ is without conjugate points, then we have 
\par \noindent
(i) Every constant sequence is an extremal sequence ; 
\par \noindent
(ii) Every extremal sequence is either injective or constant ; 
\par \noindent
(iii) For every $r \in \Z^d$, the quantity $S(x,x+r)$ does not depend on $x$.}

\medskip
{\it Proof} : Let $x \in \R^d$, and $(x_n)_{n \in \Z}$ the extremal sequence for which 
$x_0 = x_1 = x$. Using proposition 2.1 with $N=1$ and $r=0$, we may conclude that $(x_n)$ 
is a constant sequence, which proves $(i)$. Let $(x_n)_{n \in \Z}$ be extremal and not injective. 
We may assume that $x_0 = x_N$ with $N \in \N^*$. The constant sequence equal to $x_0$ is 
extremal, so corollary 1.5 tells us that $(x_n)$ is a constant sequence, which proves $(ii)$. 
For all $x \in \R^d$ and $r \in \Z^d$, we have $S(x,x+r) = \AA_1(x,x+r) = f(x)$, and according to 
lemma 2.4 this quantity does not depend on $x$, which proves $(iii)$.  
\hfill{$\square$}

\bigskip
\bigskip
{\Large 3. Some invariant Lagrangian submanifolds of $T^*\T^d$}

\bigskip
\medskip

In this section, we shall see how the translation-invariant orbits of $\Fb$ may be used 
to construct invariant Lagrangian graphs in $T^*\T^d$. We first introduce some notations. 
For every $r \in \Z^d$ and every $N \in \N^*$, we consider the following sets : 
$$\GGb_{N,r} = \{ (x,y) \in \TT, \ x \in \R^d \ {\rm and} \ \vp^N(x,y) = (x+r, y+r) \} 
\ \ {\rm and} \ \ \GGb^*_{N,r} = \LL \left ( \GGb_{N,r} \right ).$$ 

They are closely related to the extremal sequences studied in the preceding section. As a 
matter of fact, if $(x,y) \in \GGb_{N,r}$, then the extremal sequence $(x_n)_{n \in \Z}$ 
for which $x_0 = x$ and $x_1 = y$ satisfies $x_N = x_0 + r$ (and hence $x_{n+N} = x_n + r$ 
for every $n$ by proposition 2.1). Reciprocally, if $(x_n)$ is an extremal sequence 
for which $x_N = x_0 + r$, then $(x_0,x_1) \in \GGb_{N,r}$. As for $\GGb^*_{N,r}$, it contains 
all the $(x,p) \in \T^d \times (\R^d)^*$ given by theorem 1 if we fix $N$ and $r$ and let $x$ vary 
in $\R^d$.  

According to the results of the last section, there exists for every $x \in \R^d$ a unique 
$y \in \R^d$ for which $(x,y) \in \GGb_{N,r}$. This implies that $\GGb_{N,r}$ (and hence 
$\GGb^*_{N,r}$ as well) is a graph. Moreover $\GGb_{N,r}$ is clearly invariant by $\vp$,  
whereas $\GGb^*_{N,r}$ is invariant by $\Fb$. Note that as a consequence of corollary 2.5, 
$\GGb_{N,0} = \{ (x,x), \ x \in \R^d \}$.  

\medskip
Now consider $\GG^*_{N,r}$, the projection of $\GGb^*_{N,r}$ on $T^*\T^d$. It turns out that this 
set has many interesting properties : 

\medskip
{\bf Proposition 3.1} : {\it The set $\GG^*_{N,r}$ satisfies  

(i) it is a graph who is $F$-invariant ;  

(ii) $\forall \o \in \GG^*_{N,r}$, $F^N(\o) = \o$ ;

(iii) it is a Lagrangian submanifold of $T^*\T^d$.}

\medskip
{\it Proof } : Let $x \in \R^d$ and let be $p$ the unique element of $(\R^d)^*$ for which 
$(x,p) \in \GGb^*_{N,r}$. Condition $(i)$ implies that if $\Fb(x,p) = (x',p')$, 
then $\Fb(x+s,p) = (x'+s,p')$ 
for every $s \in \Z^d$. Therefore $(x+s,p) \in \GGb^*_{N,r}$ for every $s \in \Z^d$. 
They all have the same projection on $T^*\T^d$, so $\GG^*_{N,r}$ is a graph. It is $F$-invariant 
as $\GGb^*_{N,r}$ is $\Fb$-invariant. This proves $(i)$. It follows from the definitions that 
if $(x,p) \in \GGb^*_{N,r}$, then $\Fb^N(x,p) = (x+r,p)$. This readily implies property $(ii)$.

\medskip
We now prove that $\GGb^*_{N,r}$ (and therefore $\GG^*_{N,r}$ as well) is a smooth 
manifold. The main difficulty is to check that we can apply the implicit function theorem to  
$$\FF : (x,p) \in \TT^* \longmapsto \pi \left ( \Fb^N(x,p) \right ) - (x+r) \in \R^d.$$
\noindent
This will imply that the map sending $x \in \R^d$ to the unique $p \in (\R^d)^*$ for which 
$(x,p) \in \TT^*$ is smooth, and hence the smoothness of $\GGb^*_{N,r}$. 

So all we need to do is to verify that at every point in $\TT^*$, the differential 
of $\FF$ with respect to $p$ is invertible. Let $(x_0,p_0) \in \TT^*$, $(x_1,p_1) = \Fb^N(x_0,p_0)$, 
and $x_2 = \FF(x_0,p_0)$. Let  
$$i : p \in (\R^d)^* \longmapsto (x_0,p) \in \TT^*$$
\noindent 
be the canonical injection. The differential of $\FF$ with respect to $p$ at the point $(x_0,p_0)$ is  
$$ D_p\FF(x_0,p_0) : v \in T_{p_0}(\R^d)^* \longmapsto 
D\pi(x_1,p_1) \circ D\Fb^N (x_0,p_0) \circ Di(p_0) \cdot v \in T_{x_2}\R^d .$$

\noindent
Let $v$ belong to the kernel of $D_p\FF(x_0,p_0)$. Then $D\Fb^N (x_0,p_0) \circ Di(p_0) \cdot v 
\in V(x_1,p_1)$. As $Di(p_0) \cdot v \in V(x_0,p_0)$, we have $Di(p_0) \cdot v = 0$ 
(because $F$ is without conjugate points), hence $v = 0$. 

\medskip
We finally show that $\GG^*_{N,r}$ is Lagrangian. This makes use of the (positive) Green bundle, introduced  
by Bialy and McKay in [Bi-McK]. It is defined as  
$$ G(x,p) = \lim_{n \rightarrow +\infty} DF^n \left ( F^{-n}(x,p) \right ) \cdot 
 V \left ( F^{-n}(x,p) \right ) \ \subset T_{(x,p)} T^*\T^d .$$
\noindent
Each $G(x,p)$ is a Lagrangian subspace of $T_{(x,p)}T^*\T^d$. This bundle is $F$-invariant, which 
means that $DF^n G(x,p) = G \left ( F^n(x,p) \right )$ for every $n \in \Z$. Let us show that 
for every $(x,p) \in \GG^*_{N,r}$, the tangent space $T_{(x,p)}\GG^*_{N,r}$ is in fact $G(x,p)$, 
and is therefore Lagragian.

To this end, we use the following result (see [Ar]) : if $v \in T_{(x,p)}(T^*\T^d)$, then 
$$ v \notin E \Longrightarrow \lim_{n \rightarrow +\infty} 
\vv D(\pi \circ F^{-n})(x,p) \cdot v \vv = + \infty,$$
\noindent
where $\vv \cdot \vv$ is the Euclidean norm. Let $(x,p) \in \GG^*_{N,r}$ and $v \in 
T_{(x,p)}\GG^*_{N,r}$. As a consequence of $(ii)$, the restriction of $F^N$ to $\GG^*_{N,r}$ 
is the identity map, and the same is true for all $F^{-nN}$ if $n \in \Z$. Passing to the 
differential, we get $DF^{-nN}(x,p) \cdot v = v$, and hence $D(\pi \circ F^{-n})(x,p) \cdot v = 
D\pi(x,p) \cdot v$ is of constant norm. This implies that $v \in G(x,p)$. Hence 
$T_{(x,p)}\GG^*_{N,r} \subset G(x,p)$, and these two vectoriel spaces have the dimension, 
so they are equal.
\hfill{$\square$}

\bigskip
\bigskip
{\Large 4. Some results in discrete weak KAM theory}

\bigskip
\medskip

Weak KAM theory was initially developed by Mather, Mané et Fathi to study the dynamics of 
some special Hamiltonian flows. This theory was adapted to the twist maps by E. Garibaldi 
and P. Thieullen. We briefly recall the facts we shall make use of in the rest of this 
paper. We refer to [Ga-Th] for the proofs. 

\medskip  
To every generating function $S$ one can associate a real $\tS$ (called `minimizing holonomic value'). 
It is defined as  
$$\tS = {\rm Inf} \left \{ \ \liminf_{n \rightarrow +\infty} {1 \over n} S(x_0, \dots , x_n) \right \},$$
\noindent
with the infimum taken over all sequences $(x_n)_{n \in \N}$ with values in $\R^d$. One also has 
$$ \tS = {\rm Inf} _{n \geq 1} \left \{ {1 \over n} S(x_0, \dots , x_n) , \ x_0, \dots ,x_n \in \R^d 
\ {\rm with} \ x_n - x_0 \in \Z^d \right \}.$$

One usually normalizes the generating function, using $S - \tS$ instead of $S$. The action of a finite 
sequence $(x_0,\dots x_n)$ is then 
$$ \tS(x_0 , \dots x_n) = S(x_0 , \dots , x_n) - n \tS.$$
\noindent
Let us note that we now have $\tS(x_0 , \dots , x_n) \geq 0$ as soon as $x_n - x_0 \in \Z^d$, and $\tS$ 
is the smallest real number with this property. 

\medskip
The Mané potential is a function $\p : \TT \longrightarrow \R$ defined as follows : for every $(x,y) \in \TT$,  
$$ \p(x,y) =  {\rm Inf} _{n \geq 1} \left \{ \tS(x_0, \dots , x_n), \ x_0, \dots , x_n \in \R^d 
\ {\rm with} \ x_0 = x \ {\rm and} \ x_n - y \in \Z^d  \right \}.$$ 
\noindent
It is a continuous function. It is $\Z^d$-periodic with respect to each variable. It satisfies the 
triangular inequality $\p(x,z) \leq \p(x,y) + \p(y,z)$.
  
\medskip
A function $u : \R^d \longrightarrow \R$ is called a sub-action if it is $\Z^d$-periodic and if   
$$\forall x \in \R^d, \ \forall y \in \R^d, \ u(y) - u(x) \leq \p(x,y).$$
\noindent
As a consequence of the triangular inequality for $\p$, the maps $\p(x_0 , \cdot)$ et $ -\p(\cdot , x_0)$ 
are sub-actions for every $x_0 \in \R^d$.

\medskip
One can associate to $S$ a subset $\AAb$ de $\TT$ called the Aubry set : $(x,y) \in \TT$ 
belongs to $\AAb$ if for every $\e > 0$ there exists an integer $n \geq 1$ and a finite sequence 
$(x_0,x_1, \dots , x_n)$ with values in $\R^d$ for which  
$$x_n - x_0 \in \Z^d, \ \vv x - x_0 \vv < \e, \ \vv y - x_1 \vv < \e, \ {\rm and } \ 
\tS(x_0,x_1 \dots , x_n)  < \e.$$

\smallskip
The Aubry set is non-empty and closed. It is invariant by the action of $\Z^d$ : if $(x,y) \in \AAb$, 
then $(x+r,y+r) \in \AAb$ for all $r \in \Z^d$. It is also invariant by $\vp$. An important property of 
$\AAb$ is that it is a Lipschitz graph. This means that the projection on the first factor 
$\pr_1 : \AAb \longrightarrow \R^d$ is injective (hence for every $x \in \pr_1(\AAb)$, there exists 
a unique $y \in \R^d$ for which $(x,y) \in \AAb$), and that the map $x \in \pr_1(\AAb) \longmapsto 
y \in \R^d$ is Lipschitz. 

\medskip
There is a simple link between $\pr_1(\AAb)$ and the Mané potential $\p$: a point $x \in \R^d$ 
belongs to $\pr_1(\AAb)$ if and only if $\p(x,x) = 0$. If this is the case, the unique element 
$y \in \R^d$ for which $(x,y) \in \AAb$ is characterized by the relations 
$$ \p(x,y) = \tS(x,y) = S(x,y) - \tS \ \ {\rm and} \  \ \p(x,y) + \p(y,x) = 0.$$ 

\medskip
We also consider the dual Aubry set $\AAb^* = \LL(\AAb) \subset \TT^*$. It is  
a Lipschitz graph, invariant by $\Fb$. It can interpretated as the set of differentials of 
sub-actions, thanks to the following result : every sub-action $u : \R^d \longrightarrow \R$ 
is differentiable at every point $x \in \pr_1(\AAb)$, the differential being $D_xu = \LL(x,y) 
\in \AAb^*$, where $y \in \R^d$ is the unique element for which $(x,y) \in \AAb$. Finally, 
if $(x,p)\in \AAb^*$, then $(x+s,p) \in \AAb^*$ for every $s \in \Z^d$, so that we can project 
$\AAb^*$ on $T^*\T^d$ ; the result is an $F$-invariant Lispchitz graph denoted by $\AA^*$. 

\medskip
In order to construct the foliation alluded to in theorem 2, we shall consider a family of Aubry 
sets, paramerized by a cohomology class $c \in H^1(\T^d , \R)$. This is how they are defined : 
let $\o$ be a closed 1-form and $\tilde{\o}$ a lift to $\R^d$. Let us denote by 
$u : \R^d \longrightarrow \R$ a primitive of the exact 1-form $\tilde{\o}$. It is easy to check 
that the map  
$$ S_u : (x,y) \in \TT \longmapsto S(x,y) + u(x) - u(y) \in \R $$
\noindent
is a generating function. 

When we replace $S$ with $S_u$, some mathematical objects associated to $S$ will be altered, 
while others remain unchanged. For example, $S$ and $S_u$ clearly have the same extremal sequences, 
so that $\vp_u = \vp$. On the other hand, $\LL$ becomes $\LL_u = T_u^{-1} \circ \LL$, where 
$T_u$ is the translation
$$T_u : (x,p) \in \TT^* \longmapsto (x,p + Du(x)) \in \TT^*.$$
\noindent
As for $\Fb$, it is changed into $\Fb_u = T_u^{-1} \circ \Fb \circ T_u$. So if $F$ is without 
conjugate points, the same is true for $F_u$. One may check that the real $\tS_u$ only depends 
on the cohomology class $c$ of $\o$, so that it can be denoted by $\tS_c$. This gives rise to 
the $\a$-Mather function $\a : c \in H^1(\T^d,\R) \longmapsto -\tS_c \in \R$, which is both 
convex and superlinear.

As a matter of fact, the Aubry set $\AAb(S_u)$ also only depends on $c$, so it will be denoted by 
$\AAb_c$. Its dual counterpart $\LL_u(\AAb_c) = T_u^{-1}(\LL(\AA_c))$ is then invariant by 
$\Fb_u = T_u^{-1} \circ \Fb \circ T_u$. As we are more interested in 
$\Fb$-invariant subsets of $\TT^*$, it is natural to define the dual Aubry set associated 
to the cohomology class $c$ as $\AAb_c^* = \LL(\AAb_c)$. This is an $\Fb$-invariant 
Lispchitz graph. Its projection $\AA^*_c$ on $T^*\T^d$ is an $F$-invariant 
Lipschitz graph of $T^*\T^d$.

We shall make use of the following notations : if $c \in (\R^d)^*$ is a cohomology class, then 
$S_c : (x,y) \in \TT \longmapsto  S(x,y) + c(x-y) \in \R$ is its associated generating function 
and $\p_c$ the corresponding Mané potential.

\bigskip
\bigskip
{\Large 5. From periodic orbits to Aubry sets}

\bigskip
\medskip

In this section, we show that if $F$ is without conjugate points, then each of the Lagrangian 
submanifolds  $\GG_{N,r}^*$ defined in section 3 is in fact a dual Aubry set $\AA^*_c$ for a suitable 
cohomology class $c$. This is the content of the following result : 

\medskip
{\bf Proposition 5.1} : {\it Let $N \geq 1$, $r \in \Z^d$ and $u : \R^d \longrightarrow \R$ 
a smooth map such that $\GGb^*_{N,r}$ is the graph of $Du$. Then $\GG^*_{N,r} = \AA^*_c$, 
$c$ being the cohomology class of the closed 1-form induced by $Du$ on $\T^d$.}

\medskip
We first establish some special properties of the sets $\AA_c^*$ and the Mané potential $\p_c$ 
when $F$ is without conjugate points. As remarked earlier, the symplectic diffeomorphism 
$\Fb_u = T_u^{-1} \circ \Fb \circ T_u$ is then free of conjugate points as well, so that we may 
use the results obtained in sction 2, using $S_u$ instead of $S$.

\medskip
{\bf Lemma 5.2} : {\it If $F$ is without conjugate points, then $\pr_1(\AAb_c) = \R^d$ for every 
cohomology class $c$.}

\medskip
{\it Proof} : We pick $y \in \R^d$, and show that $y \in \pr_1(\AAb_c)$, i.e. $\p_c(y,y) = 0$. 
Let $x \in \pr_1(\AAb_c)$. As $\p_c(x,x) = 0$, there exists for every $\e > 0$ a finite sequenece 
$(x_0, \dots, x_n)$ with $x_0 = x$, $x_n = x_0 + r$ and $r \in \Z^d$, and 
$\tS_c(x_0, \dots , x_n) \leq \e$. We may assume that $(x_0, \dots , x_n)$ is an extremal sequence 
(see remark 1.6). Then we have (with the notations introduced in part 2) 
$S_c(x_0, \dots , x_n) = \AA_n(x,x+r) = f(x)$. Lemma 2.4 tells us that the extremal sequence 
$(y_n)_{n \in \Z}$ with $y_0 = y$ and $y_n = y_0 + r$ satisfies 
$S(y_0, \dots , y_n) = S(x_0, \dots , x_n)$. Hence   
$$S_c(y_0, \dots , y_n) = S(y_0, \dots , y_n) + c(y_0 - y_n) 
= S(x_0, \dots , x_n) + c(x_0 - x_n) =  S_c(x_0, \dots , x_n)$$
\noindent
and therefore $\tS_c(y_0, \dots , y_n) = \tS_c(x_0, \dots , x_n) \leq \e$. It follows that   
$\p_c(y,y) \leq \e$. This holds for every $\e > 0$, so that $\p_c(y,y) = 0$. \hfill{$\square$}

\medskip
{\bf Lemma 5.3} : {\it If $F$ is without conjugate points, then $\p_c$ is additive and 
antisymmetric for every cohomology class $c$: 
$$ \forall (x,y,z) \in (\R^d)^3, \  \p_c(x,z) = \p_c(x,y) + \p_c(y,z) \ \ {\rm and} \ \  
\p_c(x,y) + \p_c(y,x) = 0.$$ }
 
{\it Proof} : 
Let us fix $x$ et $y$ in $\R^d$. As explained in section 4, the maps $\p_c(x,\cdot)$ and $\p_c(y,\cdot)$ 
are two sub-actions, and are therefore differentiable at every $z \in \pr_1(\AAb_c)$, both differentials 
being equal to $\LL_c(z,z')$, with $(z,z') \in \AAb_c$. As we know that $\pr_1(\AAb_c) = \R^d$, we may 
conclude that these two maps are differentiable everywhere, with the same differential. Hence 
they are equal up to a constant :  
$$ \exists C \in \R \ {\rm s.t.} \ \forall z \in \R^d, \ \p_c(x,z) = \p_c(y,z) + C.$$
\noindent
Choosing $z = y$ and then $z = x$, we get $C = \p_c(x,y) = - \p_c(y,x)$. This yields the two relations 
$\p_c(x,y) + \p_c(y,x) = 0$ and $\p_c(x,z) = \p_c(x,y) + \p_c(y,z)$. \hfill{$\square$}

\medskip
{\bf Remark 5.4} : When $F$ is without conjugate points, the dual Aubry set $\AAb^*$ is then 
the graph of the differential of the maps $\p_c(x_0,\cdot)$, and the same is true for its 
projection on $T^*\T^d$. As for every $c \in H^1(\T^d,\R)$, we have $\AAb_c^* = \LL(\AAb_c) = 
T_u \circ \LL_u (\AAb_c)$, the set $\AA^*_c$ is the graph of a closed 1-form whose cohomology 
class is $c$.   

\medskip
We are now able to prove proposition 5.1. Let us fix $N \in \N^*$, $r \in \Z^d$, and  
$u : \R^d \longrightarrow \R$ a smooth map for which $\GGb_{N,r}^*$ is the graph of 
$Du$. The projection of $Du$ on $\T^d$ is then a closed 1-form with cohomology class $c$. 
We want to show that $\GG_{N,r}^* = \AA^*_c$.

\medskip
We first handle the case where $c = 0$, so that $u$ is a $\Z^d$-periodic function.  
For every $x \in \R^d$, $(x,Du(x)) \in \GGb_{N,r}^*$ and this set is invariant by $\Fb$, so 
we have $\Fb(x,Du(x)) = (y,Du(y))$ for a (unique) $y \in \R^d$, denoted by $y = y(x)$. 
We shall make use of the following result (see [Go] (Theorem 35.2, page 128) or [McK-Me-St] 
(Theorem 1, page 569) for a proof) : 

\medskip
{\bf Lemma 5.5} : {\it Let $u : \R^d \longrightarrow \R$ be a $C^2$ map and $\GG^* \subset \TT^*$ 
the graph of $Du$. Assume $\GG^*$ is invariant par $\Fb$ and define $\GG = \LL^{-1}(\GG^*)$. 
Then there exists a real number $C$ for which   
$$(*) \hskip 3cm \forall (x,y) \in \GG, \ S(x,y) + u(x) - u(y) = C. $$
More precisely, we have 
$$(**) \hskip 3cm \forall (x,y) \in \TT, \ u(y) - u(x) \leq S(x,y) - C,$$
\noindent
and equality holds if and only if $(x,y) \in \GG$.}

\medskip
Let $x_0 \in \R^d$, and $(x_k)$ the sequence defined by $x_{k+1}=y(x_k)$. Then $(x_k)$ 
is an extremal sequence. As $(x_0,Du(x_0)) \in \GGb^*_{N,r}$, we have $(x_0,x_1) \in \GGb_{N,r}$ 
and hence $x_N = x_0 + r$. As a consequence of $(*)$, $u(x_{k+1}) - u(x_k) = S(x_k,x_{k+1}) - C$ 
for every integer $k$. Summing up these equalities, we get  
$$\sum _{k=0}^{N-1} S(x_k,x_{k+1}) - N \times C = u(x_N) - u(x_0),$$
\noindent
and the right-hand side vanishes because $u$ is $\Z^d$-periodic. So we have 
$C = {S(x_0, \dots , x_N) \over N}$, and this implies $C \geq \tS$ by definition of $\tS$. 
Applying inequality $(**)$, we obtain  
$$\forall x \in \R^d, \ \forall y \in \R^d, \ u(y) - u(x) \leq S(x,y) - \tS = \tS(x,y),$$  
\noindent
and this means that $u$ is a sub-action. As explained in section 4, the differential of $u$ 
at every point of $\pr_1(\AAb)$ belongs to $\AAb^*$. Since $\pr_1(\AAb) = \R^d$, the graph 
of $Du$ (that is, $\GGb^*_{N,r}$) is then included in $\AAb^*$ ; as $\AAb^*$ is also a graph, 
these two sets are the same.

Assume now that $c \neq 0$. Let $S_u : (x,y) \longmapsto S(x,y) - u(x) + u(y)$ be the generating 
function and 
$$T_u : (x , p) \in \TT^* \longmapsto (x , p + Du(x)) \in \TT^*$$
\noindent
the translation. As $S$ and $S_u$ have the same extremal sequences, the sets $\GGb_{N,r}(S)$ 
and $\GGb_{N,r}(S_u)$ are equal. Using this and the fact that $\LL_u = T_u^{-1} \circ \LL$, 
we get 
$$\GGb^*_{N,r}(S_u) = \LL_u(\GGb_{N,r}(S_u)) = T_u^{-1} \circ \LL (\GGb_{N,r}(S)) 
= T_u^{-1}(\GGb^*_{N,r}(S)).$$
\noindent
The very definition of $u$ implies that this set is the null section. We may then apply the 
preceding case : the null section is in fact the dual Aubry set associated to $S_u$, and this 
means that $\GGb^*_{N,r} = \AAb^*_c$.

\bigskip
\bigskip
{\Large 6. Some supplementary results on Aubry sets}

\bigskip
\medskip

In this section, we establish some technical properties concerning Aubry sets. They will be needed 
for the proof of theorem 2. The main problem is the following : if $(c_n)$ is a sequence of 
cohomological classes that converges to $c$, what can be said about the Aubry sets $\AA^*_{c_n}$ 
and the Mané potentials $\p_{c_n}$ ? Do they converge in some sense to $\AA^*_{c}$ and $\p_{c}$ ? 
In the Hamiltonian case, every Aubry set is contained in a level set of the Hamiltonian, so 
that the $\AA^*_{c_n}$ may not explode as $n$ goes to infinity. There is no such easy argument   
in the discrete case, and therefore some new techniques are required. We state and prove four 
results ; only the last one requires $F$ to be without conjugate points.

\medskip
{\bf Lemme 6.1} : {\it Let $c$ be a cohomology class, let $(x,y) \in \AA_c$ and $(y,z) = \vp (x,y)$. 
Then} 
$$\tS_c \geq S(x,y) + S(y,z) - S(x,z).$$ 

\medskip
{\it Proof} : As $\AAb_c$ is invariant by $\vp$, both $(x,y)$ and $(y,z)$ belong to $\AAb_c$, 
so that  
$$ \p_c(x,y) = S_c(x,y) - \tS_c \ \ {\rm and} \ \  -\p_c(z,y) = \p_c(y,z) = S_c(y,z) - \tS_c. $$
\noindent
Summing up these two equalities and using the triangular inequality for $\p_c$, we get 
$$ S_c(x,y) + S_c(y,z) - 2 \tS_c = \p_c(x,y) - \p_c(z,y) \leq \p_c(x,z).$$
\noindent
By definition of $\p_c$, one has $\p_c(x,z) \leq S_c(x,z) - \tS_c$, whence    
$$ S_c(x,y) + S_c(y,z) - 2 \tS_c \leq S_c(x,z) - \tS_c,$$
\noindent
and therefore 
$$ \tS_c \geq S_c(x,y) + S_c(y,z) - S_c(x,z) = S(x,y) + S(y,z) - S(x,z).$$ 
\hfill{$\square$}

\medskip
{\bf Lemme 6.2} : {\it Let $(c_n)$ be a convergent sequence of cohomology classes, with 
$c_n \longrightarrow c$. The for every $\e > 0$, one has } 
$$\forall (x,y) \in \TT, \ \exists N_0 \in \N \ {\rm s.t.} \ 
n \geq N_0 \Longrightarrow \p_{c_n}(x,y)  \leq \e + \p_c(x,y).$$

\smallskip
{\it Proof} : Let $\e > 0$ and $(x,y) \in \TT$. By definition of $\p_c(x,y)$, there exists an 
integer $N \geq 1$ and a finite sequence $\g = (x_0, \dots , x_N)$ with $x_0 = x$, 
$y - x_N \in \Z^d$, and $\tS_c(\g) \leq \p_c(x,y) + \e$. As  
$$\tS_c(\g) = S(\g) + c \cdot (x_0 - x_N) - N \tS_c \ {\rm and} \ 
\tS_{c_n}(\g) = S(\g) + c_n \cdot (x_0 - x_N) - N \tS_{c_n},$$
\noindent
this implies, as $c_n \longrightarrow c$, that 
$$ \lim_{n \rightarrow +\infty} \tS_{c_n}(\g)  = 
\lim_{n \rightarrow +\infty} \left ( S(\g) + c_n \cdot (x_0 - x_N) - N\tS_{c_n} \right ) 
= \tS_c(\g).$$

\noindent 
So if $n$ is large enough, one has $\tS_{c_n}(\g) \leq \tS_c(\g) + \e$, 
and hence $\tS_{c_n}(\g) \leq \p_c(x,y) + 2\e$. Since  $\p_{c_n}(x,y) \leq \tS_{c_n}(\g)$, 
we finally get $\p_{c_n}(x,y) \leq \p_c(x,y) + 2\e$.  \hfill{$\square$}

\medskip
{\bf Lemma 6.3} : {\it Let $K$ be a compact set in $H^1(\T^d, \R)$. There exists a constant $M \geq 0$ 
such that }
$$ \forall c \in K, \ \forall (x,y) \in \AAb_c, \ \vv y - x \vv \leq M.$$

\smallskip
{\it Proof} : Here we use a proof by contradiction. If the conclusion was not true, we could find a 
sequence $(c_n)$ in $K$ and a sequence $(x_n,y_n)$ in $\TT$ with $(x_n,y_n) \in \AAb_{c_n}$ 
for every $n$ and $\vv y_n - x_n \vv \longrightarrow + \infty$. Since $(x_n,y_n) \in \AAb_{c_n}$, 
one has  
$$ \forall n , \ \p_{c_n}(x_n,y_n) = \tS_{c_n}(x_n,y_n) = 
S(x_n,y_n) - c_n \cdot (x_n - y_n) - \tS_{c_n}. \ \ \ \  (*)$$

\noindent
As the sequence $(c_n)$ is bounded, there exists a constant $C \geq 0$ for which  
$$ \forall n \in \N, \ \vert \tS_{c_n} \vert \leq C \ \ {\rm and} \ \ 
\vert c_n \cdot (x_n - y_n) \vert \leq C \vv x_n - y_n \vv.$$
\noindent
Since $S(x_n,y_n) \geq \a + \b \vv x_n - y_n \vv + \g \vv x_n - y_n \vv^2$ (according to lemma 1.1), 
the right hand side of $(*)$ is unbounded as $n$ goes to infinity. We will now check that the 
left hand side remains bounded, and thus get a contradiction. Since $\p_{c_n}$ is $\Z^d$-périodic 
with respect to each variable, one has  
$$ \p_{c_n}(x_n,y_n) \leq  {\rm Max} \{ \p_{c_n}(x,y), \ (x,y) \in [0,1]^d \times [0,1]^d \}. $$
\noindent
Now $\p_{c_n}(x,y) \leq \tS_{c_n}(x,y) = S(x,y) - c_n \cdot (x - y) - \tS_{c_n}$, and this quantity 
is bounded since the three variables $x$, $y$ et $c_n$ belong to compact sets. \hfill{$\square$}

\medskip
{\bf Lemma 6.4} : {\it Assume that $F$ is without conjugate points. Let $K$ be a compact set in 
$H^1(\T^d, \R)$, and let $x \in \R^d$. Then the maps $\p_c(x, \cdot)$ are uniformly Lispchitz.}

\smallskip
{\it Proof} : The maps $\p_c(x, \cdot)$ are $\Z^d$-periodic, and everywhere differentiable 
since $F$ is without conjugate points. So all we need to do is to check that the differentials 
$D_y\p_c(x,y)$ are uniformly bounded when $c \in K$ and $y \in [0,1]^d$. For every such $y$, 
let $y' \in \R^d$ with $(y,y') \in \AAb_c$. We then have $D_y\p_c(x,y) = \LL_c(y,y') = 
T_c^{-1} \circ \LL(y,y')$. To conclude the proof, simply note that $c$ is bounded, as well as 
$\vv y - y' \vv$, according lemma 6.3. \hfill{$\square$}

\bigskip
\bigskip
{\Large 7. A continuous foliation of $T^*\T^d$}

\bigskip
\medskip

In this part, we give a proof of theorem 2.  G. Paternain and M. Paternain showed in [Pa-Pa] 
that a Tonelli Hamiltonian is free of conjugate points if along every orbit of the flow there is 
a bundle of Lagrangian, flow-invariant subspaces. It is easy to adapt their proof to the twist 
maps. So if $T^*\T^d$ is foliated by Lagrangian, $F$-invariant graphs, then $F$ is without 
conjugate points. We just have to prove the converse implication.

We therefore assume that $F$ is without conjugate points and check that the dual 
Aubry sets $\AA^*_c$, with $c$ in $H^1(\T^d,\R)$, are the leaves of a continuous 
foliation of $T^*\T^d$. Let us establish that the sets $\AA^*_c$ realize a partition 
of $T^*\T^d$. We first prove that these sets are disjoint.

\medskip 
{\bf Proposition 7.1} : {\it Let $c$ and $d$ be two cohomology classes. If $c \neq d$, then  
$\AA_c^* \cap \AA_d^* = \O $.} 

\medskip 
{\it Proof} : We assume that $\AA^*_c$ and $\AA^*_d$ are not disjoint, so that $\AAb_c$ 
and $\AAb_c$ intersect at some point $(x,x') \in \TT$, and show that we then have 
$\AAb_c = \AAb_d$. This implies $\AA^*_c = \AA^*_d$ and therefore $c = d$ (see remark 5.4). 

We first prove that $\tS_{c + d \over 2} = {1 \over 2}(\tS_c + \tS_d)$. Let $(x_n)$ be the extremal 
sequence for which $x_0 = x$ and $x_1 = x'$. As the Aubry set $\AAb_c$ is invariant by $\vp$, one 
has $(x_k,x_{k+1}) \in \AAb_c$ for every integer $k$, so that 
$\p_{c}(x_k,x_{k+1}) = \tS_{c}(x_k,x_{k+1}) = S_{c}(x_k,x_{k+1}) - \tS_c$ 
for all $k$. Summing up these equalities and using the fact that $\p_c$ is additive, we get 
$$ \sum _{k=0}^{n-1}S_{c}(x_k,x_{k+1}) = n \tS_c + \p_{c}(x_0,x_n) = n \tS_c + O(1),$$
\noindent
as $\p_c$ is bounded. A similar equality holds for the cohomology class $d$, so that 
$$ \sum _{k=0}^{n-1}S_{c+d \over 2}(x_k,x_{k+1}) 
= {1 \over 2} \sum _{k=0}^{n-1} S_{c}(x_k,x_{k+1}) + S_{c}(x_k,x_{k+1}) 
= n {\tS_c + \tS_d \over 2} + O(1).$$
\noindent
This implies, by definition of $\tS_{c+d \over 2}$, that $\tS_{c+d \over 2} \leq 
{1 \over 2} (\tS_c + \tS_d)$. On the other hand, the map $c \longmapsto \tS_c$ is concave, 
hence we have equality : $\tS_{c + d \over 2} = {1 \over 2}(\tS_c + \tS_d)$.

Let us see how to use this relation to prove that $\AAb_c = \AAb_d$. Pick 
$(y,y') \in \AAb_{c+d \over 2}$. This means that for every $\e > 0$ there is a finite 
sequence $(y_0, \dots , y_n)$ with $\vv y - y_0 \vv \leq \e$, $\vv y' - y_1 \vv \leq \e$, 
$y_n - y_0 \in \Z^d$, and 
$$ \tS_{c+d \over 2}(y_0,y_1, \cdots ,y_n) = 
\sum _{k=0}^{n-1}S_{c+d \over 2}(y_k,y_{k+1}) - n \tS_{c+d \over 2} \leq \e. $$
\noindent
As $\tS_{c + d \over 2} = {1 \over 2}(\tS_c + \tS_d)$, this may be rewritten as  
${1 \over 2} ( \Sigma_c + \Sigma_d ) \leq \e$, where  
$$ \Sigma_c = \sum _{k=0}^{n-1}S_{c}(y_k,y_{k+1}) - n \tS_c  \ \ {\rm and} \ \ 
\Sigma_d = \sum _{k=0}^{n-1}S_{d}(y_k,y_{k+1}) - n \tS_d.$$
\noindent
AS $\Sigma_c$ and $\Sigma_d$ are both nonnegative quantities, each of them must be smaller than 
$2\e$. This implies that $(x,x')$ belongs to $\AAb_c$ and to $\AAb_d$. This proves that 
$\AAb_{c+d \over 2} \subset \AAb_c \cap \AAb_d$. But these three Aubry sets are all graphs 
and their projections on the first factor is $\R^d$, so they are equal. \hfill{$\square$}

\medskip
Then we establish that the union of all these dual Aubry sets $\AAb^*_c$ is equal to 
$\T^d \times (\R^d)^*$, and that they vary continuously with $c$. 

\medskip 
{\bf Proposition 7.2} : {\it For every $x \in \R^d$, the map 
$$ F_x : c \in H^1(\T^d,\R) \longmapsto  p \in (\R^d)^*, \ \ {\rm with} \ (x,p) \in \AAb^*_c,$$
\noindent
is a homeomorphism.}

\medskip 
{\it Proof} : We first establish that $F_x$ is coercive. Let $K$ be a compact set in $(\R^d)^*$, 
$c \in F_x^{-1}(K)$, $p = F_x(c)$ (so that $(x,p) \in \AAb^*_c$), $(x,x') = \LL^{-1}(x,p) \in \AAb_c$ 
and $(x',x'') = \vp((x,x')) \in \AAb_c$. According to lemma 6.1, we then have 
$$\tS_c \geq S(x,x') + S(x',x'') - S(x,x'').$$ 
\noindent 
As $p \in K$, $x'$ and $x''$ remain in compact sets in $\R^d$, so that the right-hand side is bounded 
below. Since the map $c \longmapsto -\tS_c$ is convex and superlinear, on may conclude that $c$ is 
bounded.   

\medskip
We next show that $F$ is continuous. Let $(c_n)$ be a sequence in $H^1(\T^d,\R)$. Assume that 
it converges to $c$. We have to prove that $F_x(c_n)$ goes to $F_x(c)$. Let $y_n \in \R^d$ 
with $(x,y_n) \in \AAb_{c_n}$ for every $n$. We shall establish that the sequence $(y_n)$ 
is convergent (the limit being some $y_\infty \in \R^d$) and that $(x,y_\infty) \in \AAb_c$. 
According to lemma 6.3, the sequence $(y_n)$ is bounded. So we only need to show that if 
$y_\infty$ is a cluster point of the sequence $(y_n)$, then $(x,y_\infty) \in \AAb_c$. 
This then implies that $y_\infty$ is unique (because $\AAb_c$ is a graph), and that the 
sequence converges to $y_\infty$.
 
So let us consider a convergent subsequence of $(y_n)$ (it will still be denoted by $(y_n)$ 
in order to keep notations as simple as possible), with limit $y_\infty \in \R^d$. As 
$(x,y_n) \in \AAb_{c_n}$, one has  
$$ \forall n \in \N, \ \ \p_{c_n}(x,y_n) = \tS_{c_n}(x,y_n) =  
S(x,y_n) + c_n \cdot (x-y_n) - \tS_{c_n}.$$
\noindent
When $n$ goes to infinity, the right-hand side $n$ converges to $S_c(x,y_\infty) - \tS_c$. 
The left-hand side may be rewriiten as $p_{c_n}(x,y_n) = \p_{c}(x,y_\infty)  + u_n + v_n$, 
with 
$$ u_n = \p_{c_n}(x,y_n) - \p_{c_n}(x,y_\infty) \ \ {\rm and} \ \ 
v_n = \p_{c_n}(x,y_\infty) - \p_{c}(x,y_\infty).$$

\noindent
According to lemma 6.4, the maps $\p_{c_n}(x,\cdot)$ are uniformly Lipschitz, hence the sequence 
$(u_n)$ converges to $0$. Moreover, we already know that 
$\p_{c_n}(x,y_n) \longrightarrow S_c(x,y_\infty) - \tS_c \geq \p_c(x,y_\infty)$, 
hence the sequence $(v_n)$ is convergent, its limit $\ell$ being nonnegative. On the other hand, 
lemma 6.2 tells us that for every $\e > 0$, one has 
$$ \p_{c_n}(x,y_\infty) \leq \e + \p_c(x,y_\infty), \ {\rm so} \ {\rm that} \ v_n \leq \e, $$
\noindent
when $n$ is large enough. This implies that $\ell$ has to be nonpositive, and therefore that 
$\lim v_n = 0$. So $\lim \p_{c_n}(x,y_n) = \p_c(x,y_\infty)$. From this we deduce that 
$\p_c(x,y_\infty) = S_c(x,y_\infty) - \tS_c$ and hence that $(x,y_\infty) \in \AA_c$.

\smallskip
To finish the proof, we use a topological argument : as $F_x$ is a continuous and injective 
map between two vectorial spaces of the same dimension, the invariance of domain (see [Do] 
page 567) states that $F_x$ is an open map. On the other hand, $F_x$ is a closed map, since 
it is continuous and coercive. Hence $F_x(H^1(\T^d, \R))$ is both open and closed, so it has 
to be equal to $(\R^d)^*$. Hence $F_x$ is bijective. As it is also continuous and open, it is 
a homeomorphism. \hfill{$\square$}

Another consequence of this proposition is that the map
$$\FF : (x,c) \in \T^d \times H^1(\T^d,\R) \longmapsto  F_x(c) \in T^*\T^d$$
\noindent
is continuous, and therefore the dual Aubry sets are the leaves of a continuous foliation 
of $T^*\T^d$.

\bigskip
\bigskip

\bigskip
[AABZ] ARCOSTANZO Marc, ARNAUD Marie-Claude, BOLLE Philippe, ZAVIDOVIQUE Maxime, {\it Tonelli 
Hamiltonians without conjugate points and $C^0$-integrability}, Math. Z., 280(2015), 165-194.

\bigskip
[Ar] ARNAUD Marie-Claude, {\it Green bundles and invariant submanifolds with prescribed dynamics}, 
lecture notes.

\bigskip
[Bi] BIALY Misha, {\it Rigidity  for  convex  billiards  on  the  hemisphere  and  hyperbolic  plane}, 
DCDS, 33(2013), 3903-3913.

\bigskip
[Bi-McK] BIALY M.L., MacKAY R.S., {\it Symplectic twist maps without conjugate points}, Israel Journal 
of Mathematics, 141(2004), 235-247.

\bigskip
[Bu] BUSEMANN Herbert, {\it The geometry of geodesics}, Academic Press, 1955. 

\bigskip
[Ch-Su] CHENG Jian, SUN Yisui, {\it A necessary and sufficient condition for a twist map being integrable}, 
Science in China (Series A), 39(1996), 709-717.

\bigskip
[Ga-Th] GARIBALDI Eduardo, THIEULLEN Philippe, {\it Minimizing orbits in the discrete Aubry-Mather model}, 
Nonlinearity, 24(2011), 563-611.

\bigskip
[Go] GOLE Christophe, {\it Symplectic twist maps}, World Scientific, 2001.

\bigskip
[McK-Me-St] MacKAY R.S., MEISS J.D., STARK J., {\it Converse KAM theory for symplectic twist maps}, 
Nonlinearity, 2(1989), 555-570.

\bigskip
[Pa-Pa] Paternain Gabriel, PATERNAIN Miguel, {\it On Anosov energy levels of convex Hamiltonaian 
systems}, Math. Z., 217(1994), 367-376. 

\bigskip
[St] STRUWE Michael, {\it Variational methods - Applications to Nonlinear Partial Differential Equations
and Hamiltonian Systems}, 4th edition, Springer, 2008.

\end{document}